\documentclass{amsart}

\usepackage{amsthm}
\usepackage[leqno]{amsmath}
\usepackage{latexsym,amsfonts,amssymb}
\usepackage{hyperref}

\newcommand{\numberseries}{\mdseries}   

\newlength{\thmtopspace}                
\newlength{\thmbotspace}                
\newlength{\thmheadspace}               
\newlength{\thmindent}                  

\renewcommand{\subparagraph}{\vspace*{\thmbotspace}}

\newtheoremstyle{bfupright head,slanted body}
                {\thmtopspace}{\thmbotspace}
                {\slshape}{\thmindent}{\bfseries}{.}{\thmheadspace}
                {{\numberseries \thmnumber{(#2) }}\thmnote{#3}}

\newtheoremstyle{bfupright head,upright body}
                {\thmtopspace}{\thmbotspace}
                {\upshape}{\thmindent}{\bfseries}{.}{\thmheadspace}
                {{\numberseries \thmnumber{(#2) }}\thmnote{#3}}

\newtheoremstyle{bfit head,upright body}
                {\thmtopspace}{\thmbotspace}
                {\upshape}{\thmindent}{\upshape}{.}{\thmheadspace}
                {{\numberseries\thmnumber{(#2) }}
                {\bfseries\itshape\thmnote{\negthickspace#3}}}

\newtheoremstyle{it head,upright body}
                {\thmtopspace}{\thmbotspace}
                {\upshape}{\thmindent}{\upshape}{.}{\thmheadspace}
                {{\numberseries\thmnumber{(#2) }}
                {\itshape\thmnote{\negthickspace#3}}}

\newtheoremstyle{fixed bf head,slanted body}
                {\thmtopspace}{\thmbotspace}{\slshape}
                {\thmindent}{\bfseries}{.}{\thmheadspace}
                {{\numberseries \thmnumber{(#2) }}\thmname{#1}\thmnote{ (#3)}}

\newtheoremstyle{fixed bf head,upright body}
                {\thmtopspace}{\thmbotspace}{\upshape}
                {\thmindent}{\bfseries}{.}{\thmheadspace}
                {{\numberseries \thmnumber{(#2) }}\thmname{#1}\thmnote{ (#3)}}

\newtheoremstyle{fixed bfit head,upright body}
                {\thmtopspace}{\thmbotspace}{\upshape}
                {\thmindent}{\bfseries\itshape}{.}{\thmheadspace}
                {{\numberseries \thmnumber{(#2) }}\thmname{#1}\thmnote{ (#3)}}

\newtheoremstyle{sc head,small body}
                {\thmtopspace}{\thmbotspace}
                {\small\upshape}{\thmindent}{\scshape}{.}{\thmheadspace}
                {\thmname{#1}}

\newtheoremstyle{numbered paragraph}
                {\thmtopspace}{\thmbotspace}{\upshape}
                {\thmindent}{\upshape}{}{0pt}
                {{\numberseries \thmnumber{(#2) }}}

\newtheoremstyle{unnumbered paragraph}
                {\thmtopspace}{\thmbotspace}{\upshape}
                {\parindent}{\upshape}{}{0pt}

\theoremstyle{bfupright head,slanted body}
\newtheorem{res}{}[section]             \newtheorem*{res*}{}

\theoremstyle{bfit head,upright body}
                 \newtheorem*{com*}{}

\theoremstyle{bfupright head,upright body}
\newtheorem{bfhpg}[res]{}               \newtheorem*{bfhpg*}{}

\theoremstyle{it head,upright body}
               \newtheorem*{ithpg*}{}

\theoremstyle{sc head,small body}

\theoremstyle{fixed bf head,slanted body}
\newtheorem{thm}[res]{Theorem}          \newtheorem*{thm*}{Theorem}
\newtheorem{prp}[res]{Proposition}      \newtheorem*{prp*}{Proposition}
\newtheorem{cor}[res]{Corollary}        \newtheorem*{cor*}{Corollary}
\newtheorem{lem}[res]{Lemma}            \newtheorem*{lem*}{Lemma}

\theoremstyle{fixed bf head,upright body}
       \newtheorem*{dfn*}{Definition}
     \newtheorem*{con*}{Construction}
\newtheorem{obs}[res]{Observation}      \newtheorem*{obs*}{Observation}
\newtheorem{rmk}[res]{Remark}           \newtheorem*{rmk*}{Remark}
\newtheorem{exa}[res]{Example}          \newtheorem*{exa*}{Example}
         \newtheorem*{exe*}{Exercise}
            \newtheorem{stp*}{Setup}

\theoremstyle{numbered paragraph}
\newtheorem{ipg}[res]{}

\theoremstyle{unnumbered paragraph}
\newtheorem{ipg*}{}

\newlength{\thmlistleft}        
\newlength{\thmlistright}       
\newlength{\thmlistpartopsep}   
\newlength{\thmlisttopsep}      
\newlength{\thmlistparsep}      
\newlength{\thmlistitemsep}     

\newcounter{eqc} 
  {\end{list}}%

\newcounter{prt}
\newenvironment{prt}{\begin{list}{\upshape (\alph{prt})}%
    {\usecounter{prt}%
      \setlength{\leftmargin}{\thmlistleft}%
      \setlength{\labelwidth}{\thmlistleft}%
      \setlength{\rightmargin}{\thmlistright}%
      \setlength{\partopsep}{\thmlistpartopsep}%
      \setlength{\topsep}{\thmlisttopsep}%
      \setlength{\parsep}{\thmlistparsep}%
      \setlength{\itemsep}{\thmlistitemsep}}}%
  {\end{list}}%

\newcounter{rqm}
\newenvironment{rqm}{\begin{list}{\upshape (\arabic{rqm})}%
    {\usecounter{rqm}%
      \setlength{\leftmargin}{\thmlistleft}%
      \setlength{\labelwidth}{\thmlistleft}%
      \setlength{\rightmargin}{\thmlistright}%
      \setlength{\partopsep}{\thmlistpartopsep}%
      \setlength{\topsep}{\thmlisttopsep}%
      \setlength{\parsep}{\thmlistparsep}%
      \setlength{\itemsep}{\thmlistitemsep}}}%
  {\end{list}}%

  {\end{list}}%

  {\end{list}}%

\newcounter{exercise}
  {\end{list}}%

\newenvironment{prf*}[1][Proof]{%
  \begin{proof}[\bf #1]
    \setcounter{equation}{0}
    \renewcommand{\theequation}{\arabic{equation}}}
  {\end{proof}
}

\newcommand{\pgref}[1]{(\ref{#1})}

\newcommand{\thmref}[2][Theorem~]{#1\pgref{thm:#2}}
\newcommand{\corref}[2][Corollary~]{#1\pgref{cor:#2}}
\newcommand{\prpref}[2][Proposition~]{#1\pgref{prp:#2}}
\newcommand{\lemref}[2][Lemma~]{#1\pgref{lem:#2}}

\newcommand{\exaref}[2][Example~]{#1\pgref{exa:#2}}
\newcommand{\rmkref}[2][Remark~]{#1\pgref{rmk:#2}}

\newcommand{\secref}[2][Section~]{#1\ref{sec:#2}}

\renewcommand{\eqref}[1]{\pgref{eq:#1}}

\newcommand{\thmcite}[2][?]{\cite[thm.~#1]{#2}}
\newcommand{\corcite}[2][?]{\cite[cor.~#1]{#2}}
\newcommand{\prpcite}[2][?]{\cite[prop.~#1]{#2}}
\newcommand{\lemcite}[2][?]{\cite[lem.~#1]{#2}}

\newcommand{\seccite}[2][?]{\cite[sec.~#1]{#2}}

\newcommand{\set}[1]{\{\,#1\,\}}
\newcommand{\setof}[2]{\{\,#1 \mid #2\,\}}

\newcommand{\qtext}[1]{\quad\text{#1}\quad}
 \newcommand{\qand}{\qtext{and}}
 \DeclareMathOperator*{\dsum}{\textstyle\sum}
 \newcommand{\m}{\mathfrak{m}}
 \newcommand{\is}{\cong}
\newcommand{\xra}[2][]{\xrightarrow[#1]{\;#2\;}}
\newcommand{\pows}[2][k]{#1[\mspace{-2.3mu}[#2]\mspace{-2.3mu}]}
 \newcommand{\Rm}{(R,\m)}
 \newcommand{\Rhat}{\widehat{R}}
 \renewcommand{\H}[2][]{\operatorname{H}_{#1}(#2)}
\newcommand{\bet}[3][R]{\beta^{#1}_{#2}(#3)}
 \newcommand{\Poin}[2][R]{\operatorname{P}^{#1}_{#2}(t)}
 \newcommand{\dptR}{\operatorname{depth}R}
 \newcommand{\dimR}{\operatorname{dim}R}
 \newcommand{\lgtR}{\operatorname{\ell}(R)}
 \newcommand{\edim}[1]{\operatorname{edim}#1}
 \newcommand{\codim}[1]{\operatorname{codim}#1}
 \newcommand{\Soc}[1]{\operatorname{Soc}{#1}}
 \newcommand{\SocR}{\operatorname{Soc}R}
\newcommand{\dpt}[2][R]{\operatorname{depth}_{#1}#2}
 \newcommand{\lgt}[2][R]{\operatorname{\ell}_{#1}(#2)}
\newcommand{\Hom}[3][R]{\operatorname{Hom}_{#1}(#2,#3)}
\newcommand{\Ext}[4][R]{\operatorname{Ext}_{#1}^{#2}(#3,#4)}
\newcommand{\tp}[3][R]{\nobreak{#2\otimes_{#1}#3}}

\numberwithin{equation}{res}
\newcommand{\rat}[2]{\Upsilon^{#1}_{#2}}
\newcommand{\ratmin}[1]{\Upsilon^{#1}}
\newcommand{\bas}[1]{\mu^{#1}}
\newcommand{\bass}[2][R]{\mu^{#2}(#1)}
\newcommand{\II}[1][R]{\operatorname{I}_{#1}(t)}
\newcommand{\KK}{\operatorname{K}\!}

\def\ds{\displaystyle}
\def\ts{\textstyle}
\newcommand{\e}{\varepsilon}
\newcommand{\betseq}[2][1]{\{\bet{i}{#2}\}_{i \ge #1}}
\newcommand{\bseq}[1][0]{\{b_i\}_{i \ge #1}}
\newcommand{\basseq}[1][d]{\{\bas{i}\}_{i \ge #1}}
\newcommand{\bassseq}[1][d]{\{\bass{i}\}_{i \ge #1}}

\newcommand{\xx}{\pmb{x}}
\newcommand{\yy}{\pmb{y}}
\newcommand{\x}{\pmb{x'}}
\newcommand{\h}[1]{h_{#1}}
\renewcommand{\k}{\mathsf{k}}
\newcommand{\Ek}{\operatorname{E}_R(\k)}
\newcommand{\Rmk}{(R,\m,\k)}

\newcommand{\rnk}[2][\k]{\operatorname{rank}_{#1}#2}
\newcommand{\dev}[1]{\varepsilon_{#1}}
\newcommand{\floor}[1]{\lfloor#1\rfloor}
\newcommand{\fiber}[2]{#1 \times_\k #2}
\renewcommand{\le}{\leqslant}
\renewcommand{\ge}{\geqslant}

\begin{document}

\title[Growth in the minimal injective resolution of a local
ring]{Growth in the minimal injective resolution\\ of a local ring}

\author[L.~W.~Christensen]{Lars Winther Christensen}

\address{L.~W.~Christensen, Department of Mathematics and Statistics,
  Texas Tech University, Lubbock, Texas 79409, U.S.A.}

\email{lars.w.christensen@ttu.edu}

\urladdr{http://www.math.ttu.edu/{\tiny $\sim$}lchriste}

\thanks{This work started while L.W.C.\ visited the University of
  Nebraska-Lincoln, partly supported by a grant from the Carlsberg
  Foundation. J.S.\ was supported by NSF grant DMS~0201904.}%

\author[J.~Striuli]{Janet Striuli}

\address{J.~Striuli, Department of Mathematics and Computer Science,
  Fairfield University, Fairfield, Connecticut 06824, U.S.A.}

\email{jstriuli@fairfield.mail.edu}

\author[O.~Veliche]{Oana~Veliche}

\address{O.~Veliche, Department of Mathematics, University of Utah,
  Salt Lake City, Utah 84112, U.S.A.}

\email{oveliche@math.utah.edu}

\urladdr{http://www.math.utah.edu/{\tiny $\sim$}oveliche}

\date{6 July 2009}



\begin{abstract}
  Let $R$ be a commutative noetherian local ring with residue
  field~$\k$ and assume that it is not Gorenstein. In the minimal
  injective resolution of $R$, the injective envelope $E$ of the
  residue field appears as a summand in every degree starting from the
  depth of $R$. The number of copies of $E$ in degree~$i$ equals the
  $\k$-vector space dimension of the cohomology module
  $\Ext{i}{\k}{R}$. These dimensions, known as Bass numbers, form an
  infinite sequence of invariants of $R$ about which little is
  known. We prove that it is non-decreasing and grows exponentially if
  $R$ is Golod, a non-trivial fiber product, or Teter, or if it has
  radical cube zero.
\end{abstract}


\maketitle

\section{Introduction}

\noindent In this paper $R$ is a commutative noetherian local ring
with maximal ideal $\m$.  It is a pervasive theme in local algebra
that properties of $R$ can be retrieved from homological invariants of
the residue field $\k=R/\m$. A prime example is vanishing of
cohomology with coefficients in $\k$. Indeed, $R$ is regular if and
only if $\Ext{i}{\k}{\k}=0$ for all $i\gg 0$, and $R$ is Gorenstein if
and only if $\Ext{i}{\k}{R}=0$ for all $i\gg 0$.

The cohomology groups $\Ext{i}{\k}{\k}$ behave rigidly: if $R$ is
singular (i.e.\ not regular), then they are all non-zero. In this
case, one focuses on their size, captured by the sequence of
invariants $\bet{i}{\k}=\rnk{\Ext{i}{\k}{\k}}$, called the \emph{Betti
  numbers} of $\k$. Through work pioneered by Gulliksen~\cite{THG80},
the asymptotic behavior of these Betti numbers is understood~well
enough to provide valuable information about singular local rings. The
sequence $\betseq[0]{\k}$ is known to be non-decreasing; it is
eventually constant if and only if $R$ is a hypersurface, and it has
polynomial growth if and only if $R$ is a complete intersection. If
$R$ is not a complete intersection, then the Betti numbers are
increasing and they grow exponentially.

As shown by Foxby~\cite{HBF71}, also the cohomology groups
$\Ext{i}{\k}{R}$ behave rigidly: if $R$ is not Gorenstein, then
$\Ext{i}{\k}{R}$ is non-zero for all $i \ge \dptR$. Their size is
captured by the \emph{Bass numbers} $\bass{i} = \rnk{\Ext{i}{\k}{R}}$,
but these invariants are not understood nearly as well as the Betti
numbers $\bet{i}{\k}$. It is not even known if the sequence
$\bassseq[\dptR]$ is non-decreasing. Existence of a local ring which
is not Gorenstein and whose Bass numbers grow polynomially is also an
open~question.

Let $M$ be a finitely generated $R$-module. For $i\ge 0$ the
$i^\mathrm{th}$ \emph{Betti number} of $M$ is defined as
$\bet{i}{M}=\rnk{\Ext{i}{M}{\k}}$. Let $\Rhat$ denote the $\m$-adic
completion~of~$R$. Foxby \cite{HBF77a} shows that there is a finitely
generated $\Rhat$-module $\Omega$ such that one has
\begin{equation}
  \label{eq:omega}
  \bass{\dimR+i} = \bet[\Rhat]{i}{\Omega}\quad\text{for all } i \ge 0.
\end{equation}
Thus, the non-zero Bass numbers of $R$, except the first
\mbox{$c:=\dimR - \dptR$}, can be realized as the Betti numbers of a
module. In particular, the sequence $\bassseq$ has at most exponential
growth. If $R$ is Cohen--Macaulay, then $c=0$ and $\Omega$ is the
canonical module for~$\Rhat$.  Jorgensen and Leuschke~\cite{DAJGJL07}
take this approach to prove that the sequence $\bassseq[\dptR]$ grows
exponentially and is eventually increasing for certain families of
Cohen--Macaulay rings.  Their work was motivated by a question of
Huneke about the asymptotic behavior of these numbers; however, they
are naturally lead to raise a question about the \emph{initial}
behavior \cite[2.6]{DAJGJL07}:

\begin{res}[Question]
  \label{qst:1}
  \it Let $R$ be a Cohen--Macaulay local ring of depth $d$. If the
  inequality $\bass{d+1} \le \bass{d}$ holds, is then $R$ Gorenstein?
\end{res}

\noindent There is already a history in local algebra for studying the
initial behavior of Bass numbers---in particular, the number
$\bass{\dimR}$ which is called the \emph{type} of
$R$. Bass~\cite{HBs62} characterized Gorenstein rings as
Cohen--Macaulay rings of type $1$.  Later, Foxby~\cite{HBF77a} and
Roberts~\cite{PRb83} proved that any local ring with $\bass{\dimR}=1$
is Gorenstein. This confirmed a conjecture of Vasconcelos~\cite{dtmc}.
\begin{equation*}
  * \, * \, *
\end{equation*}
\noindent In this paper we study the initial, local, and asymptotic
behavior of Bass numbers for certain families of local rings; by
\emph{local} behavior we mean comparison of a few consecutive
numbers. While the realization of Bass numbers as Betti numbers
\eqref{omega} remains important for our work, it provides no
information about the initial behavior of the Bass numbers of a ring
that is not Cohen--Macaulay. We employ a bag of ad hoc methods to deal
with this situation.

The \emph{embedding dimension} of $R$, denoted $\edim{R}$, is the
minimal number of generators of the maximal
ideal. Question~\pgref{qst:1} is answered affirmatively in
\cite{DAJGJL07} when $\edim{R}$ is $d+2$, where $d$ denotes the depth
of $R$. This is the first interesting, case inasmuch as $R$ is a
hypersurface if $\edim{R} \le d+1$.  We improve the result from
\cite{DAJGJL07} by computing the Bass numbers of these rings in closed
form; see \rmkref[]{cm codim 2}. The conclusion is strong: if $R$ is
not Gorenstein, then there are (in)equalities
\begin{equation*}
  \bass{d+1} = \bass{d} + 1 \quad\text{and}\quad \bass{i+1} \ge
  2\bass{i} \ \text{for all $i\ge d+1$}.
\end{equation*}
Similar detailed information is obtained for other rings. The Main
Theorem below extracts the overall conclusion that applies to several
families of rings. For the Cohen--Macaulay rings among them, it
answers Question~\pgref{qst:1} affirmatively.

\begin{res}[Main Theorem]
  \label{thm:main}
  Let $\Rmk$ be a local ring of embedding dimension $e$ at least $2$
  and depth $d$. Assume that $R$ satisfies one of the following
  conditions:
  \begin{rqm}
  \item \label{G} $R$ is Golod and it has $e-d> 2$ or $\bass{d} > 1$.
  \item \label{F} $R$ is the fiber product of two local rings (both
    different from $\k$) and not Golod.
  \item \label{S} $R$ is artinian with $\Soc{R}\not\subseteq \m^2$.
  \item \label{C} $R$ is not Gorenstein and $\m^3=0$.
  \item \label{T} $R$ is Teter; that is, $R \is Q/\Soc{Q}$ where $Q$
    is artinian and Gorenstein.
  \end{rqm}
  Then the sequence of Bass numbers $\bassseq$ is increasing and has
  exponential growth; if $R$ satisfies $(1)$, $(3)$, $(4)$, or $(5)$,
  then the growth is termwise exponential.
\end{res}

\noindent Notice that the assumption $\edim{R} \ge 2$ only excludes
hypersurface. A sequence $\{a_i\}_{i\ge0}$ is said to have
\emph{exponential growth} if there exists a real number $A >1$ such
that $a_i\ge\ A^i$ for all $i\gg 0$, and the growth it said to be
\emph{termwise exponential} of \emph{rate} $A$ if there exists a real
number $A>1$ such that $a_{i+1}\ge Aa_{i}$ for all $i\gg 0$.

Golod rings and fiber products are, in general, far from being
Gorenstein, so one could expect their Bass numbers grow rapidly. Teter
rings have been called ``almost Gorenstein''~\cite{CHnAVr06}, but they
still differ significantly from Gorenstein rings, and some of them are
even Golod. In this perspective, \pgref{C} is the most surprising part
of the Main Theorem, as there is empirical evidence in \cite{DAn82}
and other works that rings with \mbox{$\m^3=0$} are excellent grounds
for testing homological questions in local~algebra.

Based on the Main Theorem---and \exaref{1} below, which shows that
two consecutive Bass numbers \emph{can} be equal and non-zero---we
extend and explicitly state the question that motivated Jorgensen and
Leuschke's work \cite{DAJGJL07}:

\begin{res}[Question]
  \it Let $R$ be a local ring of depth $d$ and assume that it is not
  Gorenstein.  Does the sequence of Bass numbers $\bassseq$ then have
  exponential growth, and is it non-decreasing? Is it eventually
  increasing, and if so, from which step?
\end{res}

Our work towards the Main Theorem started from two explicit
computations of Bass numbers. We discovered \exaref[]{1} while
computing examples with the aid of \textsc{Macaulay~2}~\cite{M2}; it
can also be deduced from work of Wiebe \cite{HWb69}.

\begin{exa}
  \label{exa:1}
  Let $\mathsf{F}$ be a field. The first few Bass numbers of the local
  ring $\pows[\mathsf{F}]{x,y}/(x^2,xy)$ are $1,2,2,4,6,10$.
\end{exa}

\noindent
This ring is not Cohen--Macaulay, so the example says nothing about
Question~\pgref{qst:1}; it merely frames it. The ring is the archetype
of the exceptional case in part \pgref{G} of the Main Theorem.  We
show in \prpref{mugolod} that the Bass numbers of such rings~(Golod
with $e=d+2$ and $\bass{d}=1$) are given by the Fibonacci numbers as
follows: $\bass{d+i} = 2F_i$ for all $i\ge 1$. In particular, the
sequence $\bassseq$ is non-decreasing with termwise exponential
growth, and it increases from the third step.

Another simple example \cite[(10.8.2)]{LLAOVl07} provides a textbook
illustration of termwise exponential growth. The Main Theorem
generalizes it in several directions.

\begin{exa}
  \label{exa:2}
  For a local ring $(R,\m)$ with $\m^2=0$ and embedding dimension $e\ge
  2$, the Bass numbers are $\bass{0}=e$ and $\bass{i} =
  e^{i-1}(e^2-1)$ for all $i \ge 1$.
\end{exa}

\noindent These rings are in the intersection of the five families in
the Main Theorem. Indeed, they are Golod and they trivially satisfy
$\m^3=0$, so they belong to \pgref{G} and \pgref{C}. Since they have
$\SocR = \m$, they belong to \pgref{S} and, as will be explained in
\pgref{Peeva}, also to \pgref{F}. Finally, it is a result of
Teter~\cite{WTt74} that a local ring with $\m^2=0$ is Teter.
\begin{equation*}
  * \, * \, *
\end{equation*}
The organization of the paper follows the agenda set by the Main
Theorem. The Appendix has results on local and asymptotic behavior of
Betti numbers for modules over artinian rings. These are used in the
proofs of parts \pgref{S} and \pgref{C} of the Main Theorem, which
make heavy use of the realization of Bass numbers as Betti numbers
\eqref{omega}. A curious upshot---an immediate consequence of
\lemref{divides}---is a reformulation of Bass' characterization of
Gorenstein rings:

\begin{res}[Characterization]
  \label{dividesc}
  If $R$ is a Cohen--Macaulay local ring and $\bass{n}=1$ for some
  $n\ge0$, then $R$ is Gorenstein of dimension $n$.
\end{res}

\noindent \exaref{1} shows that the Cohen--Macaulay hypothesis is
necessary in \pgref{dividesc}, but it would be interesting to know if
there is a similar reformulation of the result of Foxby~\cite{HBF77a}
and Roberts~\cite{PRb83}. The best one can hope for is a positive
answer to:

\begin{res}[Question]
  \it Let $R$ be a local ring. If $\bass{n}=1$ for some $n\ge \dimR$,
  is then $R$ Gorenstein of dimension $n$?
\end{res}
\noindent If $R$ is an integral domain, then an affirmative answer is
already contained in \cite{PRb83}.

\section{Golod rings}
\label{sec:golod}

\noindent
The central result of this section, \thmref{mugolod}, is part
\pgref{G} of the Main Theorem.  Throughout the section, $d$ denotes
the depth of $R$ and $e$ its embedding dimension.

We use the standard notation for Koszul homology: given a sequence
$\yy$ of elements in the maximal ideal $\m$ and an $R$-module $M$, the
$i^\mathrm{th}$ homology module of the Koszul complex
$\tp{\KK{(\yy)}}{M}$ is denoted $\H[i]{\yy;M}$. Moreover, the notation
$\H[i]{\yy;R}$ is abbreviated $\H[i]{\yy}$; see also~\cite[IV.A.\S
1]{localg}.

\begin{ipg}
  Golod rings are the local rings for which the Betti numbers of the
  residue field have extremal growth; see \seccite[5]{ifr}. All
  hypersurface rings are Golod, and a Golod ring is Gorenstein if and
  only if it is a hypersurface; see \cite[rmk.~after
  prop.~5.2.5]{ifr}. The \emph{codimension} of $R$ is defined as
  $\codim{R} = \edim{R} - \dimR$. Every ring of codimension at most
  $1$ is Golod; see \prpcite[5.2.5]{ifr}. Thus, the ring in \exaref{1}
  is Golod. So is the ring in \exaref{2}, but for a different reason;
  see \prpcite[5.3.4.(1)]{ifr}.

  Here we use a characterization of Golod rings in terms of Bass
  numbers. The Bass numbers of $R$ are encoded into a formal
  power~series,
  \begin{equation*}
    \II = \sum_{i=0}^{\infty}\bass{i}t^{i},
  \end{equation*}
  called the \emph{Bass series} of $R$.

  Assume that $R$ is singular, and let $\xx$ be a minimal system of
  generators for~$\m$. Avramov and Lescot \cite[(0.2)]{LLAJLs82} prove
  that there is a coefficient-wise inequality
  \begin{equation}
    \label{eq:golodI}
    \II \preccurlyeq \frac{\sum_{i=0}^{e-1}\rnk{\H[e-i]{\xx}}
      t^i-t^{e+1}}{1 -\sum_{i=1}^e \rnk{\H[i]{\xx}} t^{i+1}},
  \end{equation}
  where equality holds if and only if $R$ is Golod.
\end{ipg}

A crucial step in the proof of \thmref{mugolod} is a reduction of
\eqref{golodI}, which comes about because $-1$ is a common root of the
numerator and the denominator. One can deduce this from the work of
Avramov, Iyengar, and Miller~\cite{AIM-06}. In \lemref{Ih} we provide
a direct argument; first we split the coefficients $\rnk{\H[j]{\xx}}$
in \eqref{golodI}:

\begin{lem}
  \label{lem:h}
  Let $R$ be singular, and let $\xx = x_1,\dots,x_e$ be a minimal
  system of generators of\, $\m$ such that $\x = x_1,\dots,x_{e-1}$
  generates an $\m$-primary ideal. For every integer $i$ the Koszul
  homology module $\H[0]{x_e;\H[i]{\x}}$ is a finite dimensional
  $\k$-vector space, and it is non-zero if and only if $i\in
  \set{0,\dots, e-d-1}$.  Moreover, for every $i \ge 0$ there is an
  equality
  \begin{equation*}
    \rnk{\H[i]{\xx}} = \rnk{\H[0]{x_e;\H[i]{\x}}} +
    \rnk{\H[0]{x_e;\H[i-1]{\x}}}.
  \end{equation*}
\end{lem}

\begin{prf*}
  Since $\x$ generates an $\m$-primary ideal, the homology module
  $\H[i]{\x}$ has finite length for every $i$, and by depth
  sensitivity and rigidity it is non-zero if and only if $0 \le i \le
  e-d-1$; see \prpcite[IV.3]{localg}.  For every $i$ and $j$ the
  module $\H[j]{x_e;\H[i]{\x}}$ is annihilated by $\m$ and hence it is
  a finite dimensional $\k$-vector space. By Nakayama's lemma,
  $\H[0]{x_e;\H[i]{\x}}$ is non-zero if and only if $\H[i]{\x}$ is so.

  For every $i$, there is a short exact sequence of finite dimensional
  $\k$-vector spaces
  \begin{equation*}
    0 \to \H[0]{x_e,\H[i]{\x}} \to \H[i]{\xx} \to
    \H[1]{x_e,\H[i-1]{\x}} \to 0;
  \end{equation*}
  see \prpcite[IV.1]{localg}. It yields
  \begin{equation*}
    \rnk{\H[i]{\xx}} = \rnk{\H[0]{x_e;\H[i]{\x}}} +
    \rnk{\H[1]{x_e;\H[i-1]{\x}}}.
  \end{equation*}
  To finish the proof we need to verify the equality
  \begin{equation*}
    \rnk{\H[1]{x_e,\H[i-1]{\x}}} =
    \rnk{\H[0]{x_e,\H[i-1]{\x}}}.
  \end{equation*}
  It follows from a length count in the exact sequence
  \begin{equation*}
    0 \to \H[1]{x_e,\H[i-1]{\x}} \to
    \H[i-1]{\x} \xra{x_e} \H[i-1]{\x} \to
    \H[0]{x_e,\H[i-1]{\x}} \to 0.\qedhere
  \end{equation*}
\end{prf*}

\begin{lem}
  \label{lem:Ih}
  Let $R$ be singular and let $\xx = x_1,\dots,x_e$ be a minimal
  system of generators of\, $\m$ such that $\x = x_1,\dots,x_{e-1}$
  generates an $\m$-primary ideal.  For $i\ge 0$ set $\h{i} =
  \rnk{\H[0]{x_e,\H[i]{\x}}}$. Then there is a coefficient-wise
  inequality
  \begin{equation*}
    \II \preccurlyeq 
    \frac{\ts\sum_{i=0}^{e-d-1}\h{e-d-1-i}t^{d+i}-t^e}{\ts
      1-\sum_{i=0}^{e-d-1} \h{i}t^{i+1}} ,
  \end{equation*}
  and equality holds if and only if $R$ is Golod.
\end{lem}

\begin{prf*}
  Since $R$ is singular, we have $e-d \ge 1$. For $i\ge 0$ set $c_i =
  \rnk{\H[i]{\xx}}$, then \eqref{golodI} takes the form
  \begin{equation*}
    \II \preccurlyeq \frac{\ts\sum_{i=0}^{e-1}c_{e-i} t^i-t^{e+1}}{\ts
      1-\sum_{i=1}^e c_it^{i+1}},
  \end{equation*}
  and equality holds if and only if $R$ is Golod.  We first verify
  that $-1$ is a root of both the numerator and denominator. Indeed,
  by \lemref{h} there are equalities $c_{e-i} = \h{e-i} + \h{e-i-1}$;
  in particular, $\h{0} = c_0 =1$. Now we have
  \begin{align*}
    \dsum_{i=0}^{e-1}c_{e-i}(-1)^i - (-1)^{e+1} & =
    \dsum_{i=d+1}^{e-1}(\h{e-i} + \h{e-i-1})(-1)^i-(-1)^{e+1}\\ &=
    \h{0}(-1)^{e-1} - (-1)^{e+1} = 0
  \end{align*}
  and
  \begin{equation*}
    1-\dsum_{i=1}^e c_i(-1)^{i+1} = 1-\dsum_{i=1}^{e-d}
    (\h{i}+\h{i-1})(-1)^{i+1} = 1 - \h{0}(-1)^2 =0.
  \end{equation*}
  Cancellation of the common factor $1+t$ gives the equality
  \begin{equation*}
    \frac{\ts\sum_{i=0}^{e-1}c_{e-i} t^i-t^{e+1}}{\ts
      1-\sum_{i=1}^e c_it^{i+1}} =
    \frac{\ts\sum_{i=0}^{e-d-1}\h{e-d-1-i}t^{d+i}-t^e}{\ts
      1-\sum_{i=0}^{e-d-1} \h{i}t^{i+1}}.\qedhere
  \end{equation*}
\end{prf*}

\begin{obs}
  \label{obs:ch}
  Let $R$ be Golod and assume it is not Gorenstein---that is, not a
  hypersurface---then one has $e-d \ge 2$; see~\cite[5.1]{ifr}. Let
  $\h{i}$ for $i\ge 0$ be as defined in \lemref{Ih}. As $\h{0}=1$ the
  Bass series of $R$ takes the form
  \begin{equation}
    \label{eq:I}
    \II = \frac{\ts\sum_{i=0}^{e-d-2}\h{e-d-1-i}t^{d+i} + t^{e-1} -  t^e}{\ts
      1 - t - \sum_{i=1}^{e-d-1} \h{i}t^{i+1}}.
  \end{equation}
  Set $\bas{i} = \bass{i}$ for $i\ge 0$. It is straightforward to
  deduce the next equalities from \eqref{I}; one can also extract them
  from the proof of \cite[(0.2)]{LLAJLs82}.
  \begin{align}
    \label{eq:mud}
    \bas{d} &= \h{e-d-1},\\
    \label{eq:mue}
    \bas{e} &= \bas{e-1} + \sum_{i=0}^{e-d-2} \bas{d+i}\h{e-d-1-i}
    -1,\quad
    \text{and} \\
    \label{eq:mue+n}
    \bas{e+n} &= \bas{e+n-1} + \sum_{i=n}^{n+e-d-2}
    \bas{d+i}\h{e-d-1+n-i} \text{ for } n\ge 1.
  \end{align}
\end{obs}

The expression for the rate of growth $A$ in the next theorem is
inspired by Peeva's proof of \prpcite[3]{IPv98}; see also
\thmcite[5.3.3.(5)]{ifr}. In view of \eqref{omega} it follows from the
latter result that the Bass sequence for a Golod ring $R$ with $e-d
\ge 2$ has termwise exponential growth. The next theorem and
\prpref{mugolod} explains the initial behavior of these Bass
sequences.

\begin{thm}
  \label{thm:mugolod}
  Let $R$ be Golod of depth $d$ and embedding dimension $e$.  If one
  of the inequalities \mbox{$e-d>2$} or $\bass{d}>1$ holds, then the
  sequence $\bassseq$ is increasing, and it has termwise exponential
  growth of rate
  \begin{equation*}
    A = \min\left\{ \frac{\bass{d+1}}{\bass{d}}, \frac{\bass{d+2}}{\bass{d+1}},
      \dots, \frac{\bass{e}}{\bass{e-1}} \right\} > 1.
  \end{equation*}
\end{thm}

\begin{rmk} In the exceptional case with $e-d=2$ and
  $\bass{d}=1\footnote{\,Let $\mathsf{F}$ be a field. The ring
    $R=\pows[\mathsf{F}]{x,y}/(x^2,xy)$ is an example. Indeed, it has
    Krull dimension~$1$ and depth~$0$, so it is Golod by
    \prpcite[5.2.5]{ifr}, and $\Hom{\k}{R}$ is generated by $x$, so
    $\bass{0} =1$.}$, the Bass numbers of $R$ are given by $\bass{d+i}
  = 2F_i$ for $i \ge 1$, where $F_i$ is the $i^\mathrm{th}$ Fibonacci
  number. In particular the sequence $\bassseq$ is non-decreasing and
  it has termwise exponential growth; see~\prpref{mugolod} and the
  remark that follows it.
\end{rmk}

\begin{prf*}[Proof of Theorem \pgref{thm:mugolod}]
  It follows from the assumptions on $R$ that it is not a
  hypersurface, so we have $e-d \ge 2$. For $i\ge 0$ set $\bas{i} =
  \bass{i}$ and adopt the notation from \lemref{Ih}. There is a
  coefficient-wise inequality
  \begin{align*}
    \ts(1-t)\II &\succcurlyeq \left(1 - t - \dsum_{i=1}^{e-d-1}
      \h{i}t^{i+1}\right)\II\\
    &= \dsum_{i=0}^{e-d-2}\h{e-d-1-i}t^{d+i} + t^{e-1}-t^e,
  \end{align*}
  where the equality follows from \eqref{I}. In particular, there are
  the following inequalities among the coefficients of $\II$:
  \begin{equation*}
    \bas{d} < \bas{d+1} < \cdots < \bas{e-1}.
  \end{equation*}
  Moreover, at least one of the inequalities $\bas{d} \ge 2$ or $e-d-
  2 \ge 1$ holds, so \eqref{mue} yields $\bas{e-1} < \bas{e}$ and,
  therefore, $A>1$. By recursion it now follows from \eqref{mue+n} and
  \eqref{mue} that $\bas{e+n} \ge A\bas{e+n-1}$ for every $n\ge 1$.
\end{prf*}

\begin{cor}
  \label{cor:Golod codim 2}
  If $R$ is Golod of codimension at least $2$, then the sequence of
  Bass numbers $\bassseq$ is increasing and has termwise exponential
  growth.
\end{cor}

\begin{prf*}
  By assumption there is an inequality $e-d \ge 2$; in particular $R$
  is not a hypersurface and hence not Gorenstein. If equality holds,
  then $R$ is Cohen--Macaulay, and then one has $\bass{d} >1$ by
  \thmcite[6.3]{HBs62}. The statement now follows from
  \thmref{mugolod}.
\end{prf*}

\begin{rmk}
  \label{rmk:cm codim 2}
  This corollary covers Cohen--Macaulay rings of codimension $2$ that
  are not Gorenstein. Indeed, such rings are Golod by \cite{GSc64};
  see also \prpcite[5.3.4]{ifr}. Also the next proposition applies to
  Cohen--Macaulay rings of codimension 2.
\end{rmk}

\begin{prp}
  \label{prp:mugolod}
  Let $R$ be of depth $d$ and embedding dimension $e=d+2$; set
  $r=\bass{d}$. If $R$ is not a complete intersection, then there is
  an equality
  \begin{equation*}
    \II = t^d\frac{r+t-t^2}{1-t-rt^2}.
  \end{equation*}
  That is, the Bass numbers of $R$ are
  \begin{equation*}
    \bass{d+i} = 
    \begin{cases}
      0 & \text{for $i < 0$}\\
      r & \text{for $i = 0$}\\
      r+1 & \text{for $i = 1$}\\
      r(r+1) & \text{for $i = 2$}\\
      2r(r+1) & \text{for $i = 3$}\\
      \bas{d+i-1} + r\bas{d+i-2} & \text{for $i \ge 4$.}
    \end{cases}
  \end{equation*}
  Set $\delta = \sqrt{1+4r}$; for $i\ge 3$ the expression for
  $\bass{d+i}$ in closed form is then
  \begin{equation*}
    \bass{d+i} =
    \frac{r+1}{2\delta}\left((2r-1+\delta)\left(\frac{1+\delta}{2}\right)^{i-1}
      + (1-2r+ \delta)\left(\frac{1-\delta}{2}\right)^{i-1}\right).
  \end{equation*}
\end{prp}

\begin{prf*}
  By \cite{GSc64} $R$ is Golod; see also \prpcite[5.3.4]{ifr}. The
  expression for the Bass series, therefore, follows from
  \eqref{I}. For $i\ge 0$ set $\bas{i} = \bass{i}$. A straightforward
  computation yields the expressions for $\bas{d+1}$ and $\bas{d+2}$,
  and \eqref{mue+n} yields the recurrence relation
  \begin{equation*}
    \bas{d+i} = \bas{d+i-1} + r\bas{d+i-2} \quad\text{for $i \ge 3$}.
  \end{equation*}
  The corresponding matrix $\binom{0\ 1}{r\ 1}$ is diagonalizable with
  eigenvalues $\frac{1}{2}(1\pm\sqrt{1+4r})$, and the expression for
  $\bas{d+i}$ in closed form follows.
\end{prf*}

\begin{rmk}
  \label{rmk:mugolod}
  Let $R$ be as in \prpref{mugolod}. If $r=\bass{d}$ is $1$, then one
  has $\bass{d+1} = 2 = \bass{d+2}$, and the recurrence relation
  yields $\bass{d+i} = 2F_i$ for $i \ge 1$, where $F_i$ is the
  $i^\mathrm{th}$ Fibonacci number.

  If $r\ge 2$, then the recurrence relation and the equality
  $\bass{d+3} = 2\bass{d+2}$ immediately yield $\bass{i+1} \ge
  2\bass{i}$ for every $i \ge d+3$, and equality holds if an only if
  $r=2$. Thus, if $r=2$, then one has $\bass{d+i} = 3(2^{i-1})$ for
  all $i \ge 1$.
\end{rmk}


\section{Fiber product rings}
\label{sec:fibers}

\noindent
In this section, $S$ and $T$ are local rings with the same residue
field $\k$ and both different from $\k$. The fiber product
$\fiber{S}{T}$ is a local ring with residue field $\k$ and embedding
dimension $e=\edim{S} + \edim{T}$. We denote its depth by $d$.

We start by observing a few fiber product rings that fail to have
increasing Bass numbers, because they are either hypersurfaces or of
the type considered in \prpref{mugolod}. The main result of the
section---\thmref{fiber}---is that they are the only (non-trivial)
ones. This will establish part \pgref{F} of the Main Theorem.

\begin{ipg}
  \label{LI}
  Let $M$ be a finitely generated $R$-module. Recall that the
  \emph{Poincar\'e series} of $M$ is the formal power series
  \begin{equation*}
    \Poin{M} = \sum_{i=0}^{\infty}\bet{i}{M}t^{i} .
  \end{equation*}
  The Poincar\'e series of $\k$ as an $\fiber{S}{T}$-module was first
  computed by Kostrikin and Shafarevich \cite{AIKIRS57}:
  \begin{equation}
    \label{eq:fiberP}
    \frac{1}{\Poin[\fiber{S}{T}]{\k}} = \frac{1}{\Poin[S]{\k}} +
    \frac{1}{\Poin[T]{\k}} -1.
  \end{equation}

  Lescot computes the quotient of the Bass series and the Poincar\'e
  series of $\fiber{S}{T}$ in \thmcite[3.1]{JLs81}. For later
  reference we record some details from this work.  If $S$ and $T$ are
  both singular, then one has
  \begin{equation}
    \label{eq:fiberI}
    \frac{\II[\fiber{S}{T}]}{\Poin[\fiber{S}{T}]{\k}} = t + \frac{\II[S]}{\Poin[S]{\k}} +
    \frac{\II[T]}{\Poin[T]{\k}}. 
  \end{equation}
  If $S$ is singular and $T$ is regular of dimension $n$, then the
  formula is
  \begin{equation}
    \label{eq:fiberII}
    \frac{\II[\fiber{S}{T}]}{\Poin[\fiber{S}{T}]{\k}} = t + \frac{\II[S]}{\Poin[S]{\k}} -
    \frac{t^{n+1}}{(1+t)^{n}}.
  \end{equation}
  If $S$ and $T$ are regular of dimension $m$ and $n$, then one has
  \begin{equation}
    \label{eq:fiberIII}
    \frac{\II[\fiber{S}{T}]}{\Poin[\fiber{S}{T}]{\k}} = t - \frac{t^{m+1}}{(1+t)^{m}} -
    \frac{t^{n+1}}{(1+t)^{n}}.
  \end{equation}
\end{ipg}

The \emph{order} of a power series $\sum_{i=0}^{\infty}v_it^i$ is the
number \mbox{$\min\setof{i \ge 0}{v_i \ne 0}$}.  Note that the order
of the Bass series $\II$ is equal to $\dptR$.

\begin{rmk}
  It follows from \eqref{fiberI}--\eqref{fiberIII} that the depth of
  $\fiber{S}{T}$ is at most~$1$, and that it is $0$ if either $S$ or
  $T$ has depth $0$. That is, one has
  \begin{equation}
    \label{eq:dpt}
    d = \min\{\dpt[]{S}, \dpt[]{T}, 1\}.
  \end{equation}
\end{rmk}

\begin{obs}
  \label{exceptions}
  Let $S$ be a $1$-dimensional regular ring.  If also $T$ is regular
  of dimension $1$, then one has $e=2$ and $d=1$ by \eqref{dpt}, so
  $\fiber{S}{T}$ is a hypersurface.
  If $T$ is either a $0$-dimensional hypersurface or a $2$-dimensional
  regular ring, then one has $e - d = 2$ by \eqref{dpt}. Moreover,
  $\fiber{S}{T}$ is Golod as both $S$ and $T$ are Golod; see
  \thmcite[4.1]{JLs83}. Finally, it follows from \eqref{fiberII} and
  \eqref{fiberIII} that $\bass[\fiber{S}{T}]{d}$ is $1$ in either
  case, and then one has $\bass[\fiber{S}{T}]{d+1} = 2 =
  \bass[\fiber{S}{T}]{d+2}$ as worked out in \rmkref{mugolod}.
\end{obs}

\begin{thm}
  \label{thm:fiber}
  Let $S$ and $T$ be local rings with common residue field $\k$ and
  assume that both $S$ and $T$ are different from $\k$. If the fiber
  product ring $\fiber{S}{T}$ is not one of the three types from {\rm
    \pgref{exceptions}}, then the sequence of Bass numbers
  $\{\bass[\fiber{S}{T}]{i}\}_{i\ge d}$ is increasing and has
  exponential growth.
\end{thm}

\noindent
Notice that part \pgref{F} of the Main Theorem, stated in the
Introduction, follows from \thmref{fiber}, as the rings in
Observation~\pgref{exceptions} are Golod.

For the proof of \thmref{fiber} we need some terminology and a
technical lemma. A power series $\sum_{i=0}^{\infty}v_it^i$ of order
$n$ is said to have non-negative (or positive) coefficients if $v_i
\ge 0$ (or $v_i > 0$) for all $i\ge n$; it has non-decreasing (or
increasing) coefficients if $v_{i+1} \ge v_i$ (or $v_{i+1} > v_i$) for
all $i\ge n$.

\begin{lem}
  \label{lem:series1}
  Let $\sum_{i=0}^{\infty}c_it^i$ be a formal power series with
  increasing coefficients and assume $c_0 > 1$. Then the following
  hold:
  \begin{prt}
  \item The coefficients of the power series
    $$V(t) = \sum_{i=0}^{\infty}v_it^i:=\frac{1}{1-t^2\sum_{i=0}^{\infty}c_it^i}$$
    grow exponentially and satisfy: $$v_0 = 1,\ \ v_1 = 0,\ \ v_2=c_0,
    \ \text{ and \ $v_{i+1}>v_i$ for all~$i\ge 2$}.$$
  \item Let $W(t) = \sum_{i=0}^{\infty}w_it^i$ be a power series of
    order $0$ and assume that the series $(1-t+t^2)W(t)$ has
    non-negative coefficients. Then the series $(1-t)V(t)W(t)$ has
    order $0$ and positive coefficients in each degree except,
    possibly, in degree $1$ where the coefficient is $w_1 - w_0 \ge
    0$.
  \end{prt}
\end{lem}

\begin{rmk}
  \label{rmk:W}
  If $R$ has positive embedding dimension, then the Poincar\'e series
  $\Poin{\k}$ satisfies the condition on $W(t)$ in
  \lemref{series1}(b). Indeed, $\Poin{\k}$ is either a power series
  with non-decreasing coefficients or the polynomial $(1+t)^n$ for
  some~\mbox{$n\ge 1$}. In the first case the claim is obvious, and in
  the second case it follows from the inequality $\binom{n}{i} \le
  \binom{n}{i-1} + \binom{n}{i+1}$, which holds for all integers $i$.
\end{rmk}

\begin{prf*}[Proof of Lemma \pgref{lem:series1}]
  (a): The equality $(1-
  t^2\sum_{i=0}^{\infty}c_it^i)(\sum_{i=0}^{\infty}v_it^i)=1$
  immediately yields $v_0=1$ and $v_1=0$, and it yields $v_2 - c_0v_0
  = 0$ whence $v_2 = c_0$. For $i \ge 3$ it yields $v_i - c_{i-2} -
  \sum_{j=0}^{i-3}c_{j}v_{i-2-j} = 0$, and it follows by recursion
  that the coefficients $v_i$ are positive. Now the desired
  inequalities
  \begin{equation*}
    v_{i+1} = c_{i-1} + \dsum_{j=0}^{i-2}c_{j}v_{i-1-j} > c_{i-2} +
    \dsum_{j=0}^{i-3}c_{j}v_{i-2-j} = v_i\quad\text{for $i\ge2$}
  \end{equation*}
  follow as the sequence $\{c_i\}_{i\ge 0}$ is increasing by
  assumption.  Finally, the expression for $v_i$ yields an inequality
  $v_{i+2} > c_0v_{i}$ for each $i\ge 2$. In particular, we have
  $v_{2j} > c_0^j$ and $v_{2j-1} > c_0^{j-1}c_1$ for $j\ge 2$. As $c_1
  > c_0 > 1$ we now have $v_i > \sqrt{c_0}^{i}$ for $i\ge 3$, so the
  sequence $\{v_i\}_{i\ge 0}$ has exponential growth.

  (b): The first equality in the computation below holds as $v_1=0$.
  \begin{align*}
    (1-t)V(t)W(t) &= (1-t)\left(1+\dsum_{i=2}^{\infty}v_it^{i}\right)W(t)\\
    &= (1-t)W(t) + t^2W(t)(1-t)\dsum_{i=2}^{\infty}v_it^{i-2}\\
    &= (1-t+t^2)W(t) +
    t^2W(t)\left((1-t)\dsum_{i=2}^{\infty}v_it^{i-2} -1 \right)
  \end{align*}
  By the assumptions on $W(t)$, the first summand in the last line of
  the display, $(1-t+t^2)W(t)$, has order $0$ and non-negative
  coefficients; the coefficient in degree $1$ is $w_1-w_0$.  It
  follows from the inequality $v_2 > 1$ that the second summand has
  order $2$, and it has positive coefficients as $v_{i+1}>v_i$ for all
  $i\ge 2$.
\end{prf*}

\begin{bfhpg*}[Proof of Theorem \pgref{thm:fiber}]
  Set $R=\fiber{S}{T}$. Assume first that $R$ is Golod; then both $S$
  and $T$ are Golod by \thmcite[4.1]{JLs83}. From \eqref{dpt} and the
  equality $e=\edim{S} + \edim{T}$, it is straightforward to verify
  that the (in)equalities $e-d \le 2$ and $\bass{d}=1$ hold exactly
  when $R$ is one of the three types of rings from
  \pgref{exceptions}. If $e-d > 2$ or $\bass{d} > 1$, then it is
  proved in \thmref{mugolod} that the sequence $\bassseq$ is
  increasing and has exponential growth.

  Assume now that $R$ is not Golod. Without loss of generality we may
  assume that $S$ is not a hypersurface;
  cf.~\thmcite[4.1]{JLs83}. Rewrite \eqref{fiberP} as
  \begin{equation*}
    \Poin{\k} =
    \frac{\Poin[S]{\k}\Poin[T]{\k}}{1-(\Poin[S]{\k}-1)(\Poin[T]{\k}-1)}.
  \end{equation*}
  Combine it with \eqref{fiberI} to obtain
  \begin{equation}
    \label{eq:poifiberI}
    \II = \frac{\II[S]\Poin[T]{\k} + \II[T]\Poin[S]{\k} +
      t\Poin[S]{\k}\Poin[T]{\k}}{1-(\Poin[S]{\k}-1)(\Poin[T]{\k}-1)}
  \end{equation}
  in case $T$ is singular. And in case $T$ is regular, combine it with
  \eqref{fiberII} to obtain
  \begin{equation}
    \label{eq:poifiberII}
    \II = \frac{\II[S](1+t)^n - t^{n+1}\Poin[S]{\k} +
      t(1+t)^n\Poin[S]{\k}}{1-(\Poin[S]{\k}-1)(\Poin[T]{\k}-1)},
  \end{equation}
  where $n \ge 1$ is the dimension of $T$.  As $S$ is not a
  hypersurface, its Poincar\'e series has increasing coefficients, so
  the power series $(1-t)(\Poin[S]{\k}-1)$ has positive coefficients,
  and hence so has $(1-t)(\Poin[S]{\k}-1)(\Poin[T]{\k}-1)$. That is,
  the series $(\Poin[S]{\k}-1)(\Poin[T]{\k}-1)$ of order $2$ has
  increasing coefficients. Moreover, the degree $2$ coefficient
  $(\edim{S})(\edim{T})$ is greater than $1$ as $S$ is not a
  hypersurface. From \lemref{series1}(a) it now follows that the
  series $1/[1-(\Poin[S]{\k}-1)(\Poin[T]{\k}-1)]$ has exponential
  growth. In \eqref{poifiberI} as well as in \eqref{poifiberII} the
  numerator has non-negative coefficients, so in either case the
  sequence $\bassseq$ has exponential growth.

  If $T$ is singular, then the series $(1-t)\Poin[T]{\k}$ has
  non-negative coefficients, and as above $(1-t)\Poin[S]{\k}$ has
  positive coefficients. Thus, the numerator in the expression
  \begin{equation*}
    (1-t)\II = \frac{(1-t)(\II[S]\Poin[T]{\k} + \II[T]\Poin[S]{\k} +
      t\Poin[S]{\k}\Poin[T]{\k})}{1-(\Poin[S]{\k}-1)(\Poin[T]{\k}-1)},
  \end{equation*}
  which is derived from \eqref{poifiberI}, is a power series of order
  $d\le 1$ with positive coefficients. As the power series
  $1/[1-(\Poin[S]{\k}-1)(\Poin[T]{\k}-1)]$ has non-negative
  coefficients and order $0$, it follows that $(1-t)\II$ has positive
  coefficients in all degrees $i\ge d$. That is, the sequence
  $\bassseq$ is increasing.

  Finally, assume that $T$ is regular of dimension $n$. The Poincar\'e
  series $(1+t)^n$ and $\Poin[S]{\k}$ satisfy the condition on $W(t)$
  in \lemref{series1}(b); see also \rmkref{W}. By \eqref{poifiberII}
  the series $(1-t)\II$ of order $d\le 1$ can be expressed as a sum
  \begin{align*}
    (1-t)\II &=
    (1-t)\frac{(1+t)^n}{1-(\Poin[S]{\k}-1)(\Poin[T]{\k}-1)}\II[S]\\ &+
    (1-t)\frac{\Poin[S]{\k}}{1-(\Poin[S]{\k}-1)(\Poin[T]{\k}-1)}(t(1+t)^n
    - t^{n+1}).
  \end{align*}
  The first summand has order $\dpt[]{S} \ge d$, and by
  \lemref{series1}(b) it has non-negative coefficients. Similarly, the
  second summand has order $1$ and positive coefficients, also in
  degree $2$ as $\bet[S]{1}{\k} - \bet[S]{0}{\k} = \edim{S} - 1 >
  0$. It follows that the power series $(1-t)\II$ has positive
  coefficients in all degrees $i\ge d$, therefore, the sequence
  $\bassseq$ is increasing. \hfill\qed
\end{bfhpg*}

\begin{rmk}
  \label{Peeva}
  Let $R$ be artinian, and assume that $\k$ splits out of the maximal
  ideal $\m$---that is, $\Soc{R}$ contains a minimal generator of
  $\m$. If $R$ is not a hypersurface, i.e.\ $\edim{R}$ is at least
  $2$, then it is a non-trivial fiber product of artinian local
  rings. Indeed, let $x,y_1,\dots,y_{m}$ be a minimal set of
  generators of $\m$, such that $x$ is in $\SocR$, then there is an
  isomorphism $R \is \fiber{R/(x)}{R/(y_1,\dots,y_{m})}$. By
  \thmref{fiber} the sequence $\bassseq[0]$ is increasing and has
  exponential growth.
\end{rmk}

One can, however, do better for this particular kind of fiber
products. Let $R$ be as in \pgref{Peeva}. If $\m^2=0$ the sequence
$\bassseq[0]$ has even termwise exponential growth; see \exaref{2}. It
the next section---see \prpref{socgrowth}---the same conclusion is
reached in the case $\m^2\ne 0$, and this establishes part \pgref{S}
of the Main~Theorem.


\section{Artinian rings}
\label{sec:art}

\noindent
In this section $\Rmk$ is artinian; the injective envelope of $\k$ is
denoted~$\Ek$.  The results in this section prepare the grounds for
the proof of part \pgref{C} of the Main Theorem from the Introduction;
it is given in the next section. The last result of this section
establishes  the Main~Theorem's part \pgref{S}.

\begin{ipg}
  \label{art notation}
  Let $M$ be a finitely generated $R$-module. We write $\lgt{M}$ for
  the length of $M$ and $\Soc{M}$ for its socle. It is straightforward
  to verify the following:
  \begin{equation}
    \label{lescot}
    \text{$\k$ is a direct summand of $M$ if and only if
      $\Soc M\not\subseteq\m M$.} 
  \end{equation}
  For $i\ge 1$ we denote the $i^\mathrm{th}$ syzygy in a minimal free
  resolution of $M$ by $M_i$; we set $M_0 = M$. Since the
  differentials in a minimal free resolution are given by matrices
  with entries in $\m$, there are equalities
  \begin{equation}
    \label{soc}
    \Soc{M_{i+1}} = \SocR^{\bet{i}{M}} \quad\text{for all } i\ge 1.
  \end{equation}
\end{ipg}

\begin{ipg}
  \label{art}
  Recall the following special case of \eqref{omega}:
  \begin{equation}
    \label{eq:mb}
    \bass{i} = \bet{i}{\Ek} \quad\text{ for all $i\ge 0$}.
  \end{equation}
  If $R$ is not Gorenstein, then the equality $\lgt{\Ek} = \lgtR$ and
  \lemref{divides} yield
  \begin{equation}
    \label{eq:mu01}
    \bass{1} \ge \bass{0}.
  \end{equation}
  Recall also that one has $\rnk{\SocR}=\bass{0}$, and denote this
  number by $r$.
\end{ipg}

The next result contains a special case of \thmref{m3}, namely the one
where the socle rank $r$ exceeds the embedding dimension.

\begin{prp}
  \label{php:hmu}
  Set $h = \max\setof{i}{\m^i\not=0}$ and assume that $R$ satisfies
  $h\ge 2$. If $r > \lgt{R/\m^h}-1$, then the sequence of Bass numbers
  $\bassseq[0]$ is increasing and has termwise exponential growth of
  rate at least
  $$\frac{r}{\lgt{R/\m^h}-1}.$$
\end{prp}

\begin{prf*}
  The assumptions on $R$ force an inequality $r>2$; in particular $R$
  is not Gorenstein. Indeed, let $e$ be the embedding dimension of
  $R$, then one has $r > e$, as $h \ge 2$ by assumption. If $e$ were
  $1$, then $R$ would be a hypersurface, so also $r$ would be $1$,
  which is impossible. Hence, one has $e \ge 2$ and $r > 2$.
 
  By \prpref{s} we need only prove the inequality $\bass{1} >
  \bass{0}$, and to this end it suffices, by \eqref{mu01}, to show
  $\bass{1} \ne r$.  Assume, towards a contradiction, that one has
  $\bass{1}=r$. Set $E=\Ek$ and consider the exact~sequence $$0 \to
  E_2 \to R^{r} \to R^r \to E \to 0,$$ that comes from the minimal
  free resolution of $E$. By \pgref{soc} and additivity of length, the
  sequence yields (in)equalities
  \begin{equation*}
    r^2 = \rnk{\SocR}^r = \rnk{\Soc{E_2}} \le \lgt{E_2}=\lgt{E}.
  \end{equation*} 
  By additivity of length, the assumptions on $R$, and the containment
  $\m^h \subseteq \SocR$ there are (in)equalities
  \begin{equation*}
    \lgtR = \lgt{R/\m^h} + \lgt{\m^h} < 1 + r + \rnk{\m^h} \le 1 + 2r.
  \end{equation*}
  By the equality $\lgt{E} = \lgtR$, the last two displays combine to
  yield $r^2 \le 1+2r$, which implies $r \le 2$, a contradiction.
\end{prf*}

\begin{lem}
  \label{lem:growth}
  If\, $\k$ is a direct summand of the first syzygy of\, $\Ek$, then
  one has $\bass{1} > \bass{0}$.
\end{lem}

\begin{prf*}
  By \eqref{mu01} it is enough to rule out the possibility of an
  equality \mbox{$\bass{1}=r$}. Suppose this equality holds. For
  brevity, set $E=\Ek$.  By hypothesis, there is an isomorphism $E_1
  \is \k\oplus N$ for some finitely generated $R$-module $N$ and,
  therefore, there is an isomorphism of syzygies $E_2 \is \m\oplus
  N_1$.  From the exact sequence \mbox{$0\to E_2\to R^r\to R^r\to E\to
    0$} one obtains $\lgt\m + \lgt{N_1} = \lgt{E}$ by additivity of
  length. Since $\lgt{E} = \lgt{R} = \lgt{\m} + 1$, it follows that
  $N_1$ is a $\k$-vector space of rank $1$. Hence there is an exact
  sequence $0\to \k\to R^{r-1}\to N \to 0$, where the surjective
  homomorphism is given by a matrix with entries in $\m$, as $N$ is a
  submodule of $\m R^r$, cf.~\pgref{soc}. Thus, we have $\SocR^{r-1} =
  \k$, and the ensuing equality of ranks $r(r-1) = 1$ is absurd.
\end{prf*}

\begin{lem}
  \label{lem:SocE1}
  Assume that $R$ is not Gorenstein, and set $E=\Ek$. Let $n$ be an
  integer such that $\Soc{E_1}$ is contained in $\m^n{E_1}$. Then one
  has $n < \max\setof{i}{\m^i\not=0}$ and $\SocR \subseteq \m^{n+1}$.
\end{lem}

\begin{prf*}
  Set $h=\max\setof{i}{\m^i\not=0}$, then $\m^{h}E_1=0$ as $E_1$ is a
  submodule of $\m R^r$; and the inequality $n < h$ follows.

  Let $e_1,\dots,e_r$ be the standard basis for $R^r$. Choose a
  minimal set $\e_1,\dots,\e_r$ of generators of $E$ and consider the
  short exact sequence $0 \to E_1 \to R^r \to E \to 0$, where the
  surjection maps $e_i$ to $\e_i$.  Suppose that $\SocR$ is not
  contained in $\m^{n+1}$ and choose an element $x$ in $(\SocR)
  \setminus \m^{n+1}$. Since $E$ is a faithful $R$-module, there is an
  index $i$ such that the element $x\e_i$ is non-zero; it clearly
  belongs to $\Soc{E}$. Similarly, by the definition of $h$, it
  follows that $\m^h E$ is a non-zero submodule of $\Soc{E}$. Since
  $\rnk{\Soc{E}}=1$, there exists an element $y$ in $\m^h$ and an
  index $j$ such that $y\e_j = x\e_i$.  The element $z = xe_i - ye_j$
  is in the socle of $E_1$, as $x$ and $y$ are in the socle of $R$. As
  $y$ is in $\m^{n+1}$ but $x$ is not, the element $z$ is not in
  $\m^nE_1 \subseteq \m^{n+1}R^r$, which contradicts the assumption
  that $\Soc{E_1}$ is contained in $\m^n{E_1}$.
\end{prf*}

\begin{rmk}
  \label{rmk:kE1}
  Set $E=\Ek$. If $\k$ is a direct summand of $\m$, that is $\SocR
  \not\subseteq \m^2$, then it follows from \lemref{SocE1} that
  $\Soc{E_1} \not \subseteq \m E_1$. Thus, $\k$ is a direct summand of
  $E_1$, cf.~\pgref{lescot}, and then \lemref{growth} yields the
  inequality $\bass{1} > \bass{0}$.
\end{rmk}

Finally, \prpref{growth} together with \rmkref{kE1} gives us:

\begin{prp}
  \label{prp:socgrowth}
  Let $R$ be an artinian local ring of embedding dimension at least
  $2$. If\, $\SocR \not\subseteq \m^2$, then the sequence of Bass
  numbers $\bassseq[0]$ is increasing and it has termwise exponential
  growth.\qed
\end{prp}

\section{Rings with $\m^3=0$}
\label{sec:m3}

\noindent The theorem below is part \pgref{C} of the Main Theorem from
the Introduction.
\begin{thm}
  \label{thm:m3}
  Let $\Rm$ be local with $\m^3=0$. If $R$ is not Gorenstein, then the
  Bass sequence $\bassseq[0]$ is increasing and has termwise
  exponential growth.
\end{thm}

\begin{prf*}
  Fix the following notation
  \begin{equation*}
    a = \rnk{\m^2}, \quad e=\rnk{\m/\m^2}, \qand r=\rnk{\SocR}.
  \end{equation*}
  The embedding dimension $e$ is at least $2$, as $R$ is not a
  hypersurface, and for $e=2$ the statement is contained in
  \corref{Golod codim 2}, cf.~\rmkref{cm codim 2}. In the rest of the
  proof, we assume $e\ge 3$.  In view of \exaref{2} we can assume that
  $\m^2$ is not $0$, i.e.\ $a>0$.  Set $E=\Ek$. For $i\ge 0$ set
  $\bas{i}=\bass{i}$ and recall from \eqref{mb} that $\bas{i} =
  \bet{i}{E}$. The containment $\m^2 \subseteq \SocR$ implies an
  inequality $a \le r$.

  If $a < r$, then $\SocR$ is not contained in $\m^2$ and \rmkref{kE1}
  gives the inequality $\bas{1} > \bas{0}$. Then
  \thmref{increasing}(a) applies to the module $E$, so the sequence
  $\basseq[0]$ is increasing and it has termwise exponential growth.

  In the remainder of the proof assume $a=r$. From the equalities
  $\lgt{E} = \lgtR$ and $1 = \rnk{\Soc{E}} = \rnk{\m^2E}$ one gets
  $\rnk{\m E/\m^2E} = \lgtR - r -1 = e$. By \lemcite[3.3]{JLs85}---see
  also \eqref{betti1}---there is an inequality
  \begin{equation}
    \label{eq:bas1}
    \bas{1} \ge e\bas{0} - \rnk{\m E/\m^2E} = e(r-1),
  \end{equation}
  and equality holds if and only if $\k$ is not a direct summand of
  $E_1$.  Since \mbox{$e\ge 3$} and $r \ge 2$, an equality $\bas{1} =
  \bas{0}$ would imply $r \ge e(r-1) \ge r +2r-3 \ge r+1$, which is
  absurd. Thus, the inequality $\bas{1} > \bas{0}$ holds,
  cf.~\eqref{mu01}.  Now \thmref{increasing}(a) applies to the module
  $E$, so the sequence $\basseq[0]$ is increasing with termwise
  exponential growth except, possibly, when $a=e=r$.

  Assume now that all three invariants $a$, $e$, and $r$ are
  equal. Note that \eqref{bas1} yields
  \begin{equation}
    \label{eq:bas14}
    \bas{1} \ge (r-1)\bas{0}.
  \end{equation}
  By \lemcite[3.3]{JLs85}---see also \eqref{betti3}---there are
  inequalities
  \begin{equation}
    \label{eq:bas3}
    \bas{i+1} \ge r(\bas{i} - \bas{i-1})\quad \text{for $i \ge 1$}
  \end{equation} 
  and equality holds if $i\ge 2$ and $\k$ is not a summand of either
  syzygy $E_{i}$ and $E_{i-1}$.

  Consider the case where the common value of $a$, $e$, and $r$ is at
  least $4$. The quantity $A=\frac{1}{2}(r+\sqrt{r^2-4r})$ is then a
  real number greater than or equal to $2$. We claim that the
  inequality $\bas{i+1} \ge A\bas{i}$ holds for all $i\ge 0$. The
  proof is by induction on $i$. The base case $i=0$ is furnished by
  \eqref{bas14} and the induction step follows from~\eqref{bas3}:
  \begin{equation*}
    \bas{i+1} \ge r(\bas{i} - \bas{i-1}) \ge r(\bas{i} -A^{-1}\bas{i})
    = r(1-A^{-1})\bas{i} = A\bas{i},
  \end{equation*}
  where the last equality follows as $A$ is a solution to the equation
  $r(1-x^{-1}) = x$.

  As $e \ge 3$ we are left with only one case to consider, namely
  $a=e=r=3$.  First we prove that the sequence $\basseq[0]$ is
  increasing and that $\k$ is a direct summand of one of the first
  four syzygies of $E$. If $\k$ is a direct summand of $E_1$ or $E_2$,
  then the sequence $\basseq[0]$ is increasing by
  \thmref{increasing}(b). Assume now that $\k$ is not a direct summand
  of $E_1$, then \eqref{bas1} yields
  \begin{equation*}
    \bas{1} = e(r - 1) = 6.
  \end{equation*}
  Assume also that $\k$ is not a direct summand of $E_2$. In the
  computation $$\rnk{\m E_2} = \rnk{\Soc{E_2}} = \rnk{\SocR^6} = 18$$
  the first equality follows from \pgref{lescot} and the second from
  \pgref{soc}. Combined with length computations in the minimal free
  resolution of $E$, this gives
  \begin{align*}
    \bas{2} = \lgt{E_2} - \lgt{\m E_2} = 4\lgt{R} - 18 = 10.
  \end{align*}
  In particular, one has $\bas{2} > \bas{1}$. If $\k$ is a direct
  summand of $E_3$, then the sequence $\basseq[0]$ is increasing by
  \thmref[]{increasing}(b). If $\k$ is not a direct summand of $E_3$,
  then \eqref{bas3}~yields
  \begin{align*}
    \bas{3} = 3(\bas{2} - \bas{1}) = 12.
  \end{align*}
  If $\k$ were not a direct summand of $E_4$, then \eqref{bas3} would
  yield $\bas{4} = 3(\bas{3} - \bas{2}) = 6 < \bas{3}$ which is
  impossible by \thmref{increasing}(b). Thus, $\k$ is a direct summand
  of $E_4$, and the sequence $\basseq[0]$ is increasing by
  \thmref[]{increasing}(b) as one has $\bas{3} > \bas{2}$.

  Finally, we can conclude that the growth of the series $\basseq[0]$
  is termwise exponential. Since $\k$ is a direct summand of a syzygy
  of $E$, the radius of convergence of the power series $\Poin{E} =
  \II$ is bounded above by that of $\Poin{\k}$. The opposite
  inequality always holds by \prpcite[1.1]{JLs85}, so the two power
  series have the same radius of convergence $\rho$. As $R$ is not
  complete intersection, the sequence $\betseq[0]{\k}$ has exponential
  growth, so one has $\rho < 1$. From work of Sun~\thmcite[1.2]{LCS94}
  it now follows that the sequence $\basseq[0]$ has termwise
  exponential growth.
\end{prf*}

\begin{bfhpg}[Scholium]
  \renewcommand{\arraystretch}{1.5}%
  Let $R$ be as in \thmref{m3}. The table below gives a lower bound
  for the rate $A$ of the termwise exponential growth of the Bass
  sequence $\bassseq[0]$ in terms of the invariants $a = \rnk{\m^2}$,
  $e=\rnk{\m/\m^2}$, and $r=\rnk{\SocR}$.

  {\small
    \begin{equation*}
      \begin{array}{c|c|c|c|c|c|c|c}
        a &a=0&a<r& \multicolumn{5}{c}{a=r} \\ \cline{4-8}   
        e &   &   & e<r & \multicolumn{2}{c|}{e=r}  & e=r+1 & e \ge r+2 \\ \cline{5-6}
        r &   &   &     & 2 \le r \le 3        &  4 \le r &       &  \\
        \hline
        A &e=r& r-a+\ds\frac{r-a}{a+e} & \ds\frac{r}{e}  & 2 &
        \ds\frac{r+\sqrt{r^2-4r}^{\phantom{1}}}{2} & \ds\frac{r^2-1}{r}
        & e-r\\
      \end{array}
    \end{equation*}
  } \renewcommand{\arraystretch}{1}

  \noindent For rings with $a=e=r=3$, Backelin and
  Fr{\"o}berg~\cite{JBcRFr85} give closed form expressions for the
  possible Poincar\'e series $\Poin{\k}$, and one can verify directly
  that the radius of convergence is less than $\frac{1}{2}$.  Sun's
  \thmcite[1.2]{LCS94} then yields the lower bound $2$ for~$A$. In all
  other cases, the bound follows by inspection of \exaref{1}, the
  proof of \thmref{increasing}(a), the proof of \thmref{m3}, and
  Peeva's proof of \prpcite[3]{IPv98}.
\end{bfhpg}





\section{Teter rings}
\label{sec:teter}

\noindent
Following Huneke and Vraciu~\cite{CHnAVr06}, we say that $R$ is
\emph{Teter} if there exists an artinian Gorenstein local ring $Q$
such that $R \is Q/\Soc Q$.

\begin{ipg}
  Let $Q$ be an artinian Gorenstein local ring which is not a
  field. Denote by $\k$ its residue field. Set $R=Q/\Soc{Q}$.  The
  maximal ideal of $Q$ is isomorphic to $\Hom[Q]{R}{Q}$, which is the
  injective hull of $\k$ as an $R$-module. Thus, $R$ is Gorenstein if
  and only if $\edim{Q} = 1$, in which case both $Q$ and $R$ are
  hypersurfaces.

  If $Q$ is not a hypersurface, i.e.\ $\edim{Q} \ge 2$, then one has
  $\edim{R} = \edim{Q}$ and works of Avramov and Levin
  \thmcite[2.9]{GLLLLA78} and Herzog and Steurich
  \prpcite[1]{JHrMSt79} provide the following expression for the Bass
  series of~$R$:
  \begin{equation}
    \label{eq:Teter1}
    \II = \frac{\Poin[Q]{\k}-1}{t(1-t^2\Poin[Q]{\k})}.
  \end{equation}
  Notice the equality $\bass{0} = \edim{R}$.
\end{ipg}

\begin{ipg}
  For each positive integer $e$ set
  \begin{align*}
    \rat{e}{i} &= \binom{e}{i}^{-1}\left(\sum_{j=0}^{i+1}
      \binom{e-1}{j}\right)\qquad\text{for } i\in\{0,\dots,e\} \quad\text{and}\\
    \ratmin{e} &= \min\{\rat{e}{i} \mid 0 \le i \le e \}.
  \end{align*}

  The quantities $\ratmin{e}$ are used below to provide lower bounds
  for the growth rate of Bass numbers, and we make a
  few observations on how to compute or estimate their values. Notice
  that one has $\rat{e}{0} = e$ and $\rat{e}{e}=2^{e-1}$ for every
  $e\ge 1$. Moreover, one has $\rat{2}{1} = 1$, and for $e\ge 3$ it is
  straightforward to establish the inequalities
  \begin{align*}
    && \rat{e}{1} &\le \frac{e}{2}\\
    &\mspace{80mu}& 1 < 1 + \binom{e}{i}^{-1}\binom{e-1}{i+1} &\le
    \rat{e}{i}& &\text{for}& 1 &\le i \le
    \floor{\textstyle\frac{e}{2}} \quad\text{and} &\mspace{80mu}&\\
    && \rat{e}{i} &\le \rat{e}{i+1}& &\text{for}&
    \mspace{-12mu}\floor{\textstyle\frac{e}{2}} &\le i \le e-1.
  \end{align*}
  In particular, one has
  \begin{equation}
    \label{eq:ratmin}
    \ratmin{1} = 1 = \ratmin{2} \quad\text{and}\quad 1 < \ratmin{e} =
    \min\{\rat{e}{i}
    \mid 1 \le i \le
    \floor{\textstyle\frac{e}{2}}\} \ \text{ for } e\ge 3. 
  \end{equation}
  Direct computations yield $\ratmin{e} =
  \rat{e}{\floor{\frac{e}{2}}}$ for $e\le 9$ but $\ratmin{10} =
  \rat{10}{4}$. To find $\ratmin{e}$ for larger values of $e$ it is
  useful to know that $\rat{e}{i}$ is an upwards convex function of
  $i$; Roger W.~Barnard proved this upon request. Thus, if $n$ is an
  integer between $1$ and $\floor{\textstyle\frac{e}{2}}$ and there
  are inequalities $\rat{e}{n-1} \ge \rat{e}{n} \le \rat{e}{n+1}$,
  then one has $\ratmin{e} = \rat{e}{n}$.
\end{ipg}

\begin{lem}
  \label{lem:A}
  Let $R$ be artinian of embedding dimension \mbox{$e \ge 1$}. For
  every real number $A < \ratmin{e}$ the formal power series
  $(1-At+At^3)\Poin{\k}$ has positive coefficients, and for
  $A=\ratmin{e}$ it has non-negative coefficients.
\end{lem}

\begin{prf*}
  By \cite[7.1]{ifr} there exist non-negative integers $\dev{i}$ such
  that
  \begin{equation*}
    \Poin{\k} = \frac{\prod_{i=1}^{i=\infty}
      (1+t^{2i-1})^{\dev{2i-1}}}{\prod_{i=1}^{i=\infty}
      (1-t^{2i})^{\dev{2i}}}
    = \frac{(1+t)^{\dev{1}}}{1-t^2}\cdot\frac{\prod_{i=2}^{i=\infty}
      (1+t^{2i-1})^{\dev{2i-1}}}{(1-t^2)^{\dev{2}-1}\prod_{i=2}^{i=\infty}
      (1-t^{2i})^{\dev{2i}}},
  \end{equation*}
  and one has $\dev{2}-1 \ge 0$ as $R$ is not regular; see
  \thmcite[7.3.2]{ifr}.  By \corcite[7.1.5]{ifr} there is an equality
  $\dev{1} = e$, so the factor $F=\frac{(1+t)^{\dev{1}}}{1-t^2}$ can
  be rewritten as
  \begin{equation*}
    F = \frac{(1+t)^{e-1}}{1-t} = (1+t)^{e-1}\sum_{i=0}^{\infty} t^i =
    \sum_{i=0}^{e-2} \left(\sum_{j=0}^i \binom{{e-1}}{j} \right)t^i +
    \sum_{i={e-1}}^\infty 2^{e-1}t^i.
  \end{equation*}
  The power series $(1-At+At^3)\Poin{\k}$ has positive/non-negative
  coefficients if the series $(1-At+At^3)F$ has positive/non-negative
  coefficients. From the computation

  \begin{align*}
    (1-At+At^3)F &= \sum_{i=0}^{e+1} \left(\sum_{j=0}^i
      \binom{{e-1}}{j} - A\left(\!\binom{{e-1}}{i-2} +
        \binom{{e-1}}{i-1}\!\right)\!\right)t^i +
    \sum_{i=e+2}^\infty 2^{e-1}t^i\\
    &= 1 + \sum_{i=1}^{e+1} \left(\sum_{j=0}^i \binom{e-1}{j} -
      A\binom{e}{i-1}\!\right)t^i + \sum_{i=e+2}^\infty 2^{e-1}t^i\\
    &= 1 + \sum_{i=0}^{e} \left(\sum_{j=0}^{i+1} \binom{e-1}{j} -
      A\binom{e}{i}\!\right)t^{i+1} + \sum_{i=e+2}^\infty 2^{e-1}t^i
  \end{align*}
  it follows that $(1-At+At^3)F$ has non-negative coefficients if and
  only if $A \le \ratmin{e}$ and positive coefficients if $A <
  \ratmin{e}$.
\end{prf*}

The next result establishes part \pgref{T} of the Main Theorem.

\begin{thm}
  Let $R$ be Teter of embedding dimension $e$ at least $2$.
  \begin{prt}
  \item If $e = 2$, then $\bass{0} = 2$ and $\bass{i} = 3(2^{i-1})$
    for $i\ge 1$.
  \item If $e > 2$, then $\ratmin{e} > 1$ and the formal power series
    $(1-\ratmin{e}t)\II$ has non-negative coefficients.
  \end{prt}
  In particular, the sequence of Bass numbers $\bassseq[0]$ is
  increasing and has termwise exponential growth.
\end{thm}

\begin{prf*}
  If $e = 2$, then $R$ is Golod; see \prpcite[5.3.4]{ifr}. Since $R$
  is Teter, it follows from \eqref{Teter1} that $\bass{0} = e =
  2$. Thus, part (a) follows from \rmkref{mugolod}.

  Assume now that $e>2$, then $\ratmin{e} > 1$ by \eqref{ratmin}. To
  prove part (b), it is sufficient to show that the series
  $(1-\ratmin{e}t)(1+t\II)$ has non-negative coefficients.  To this
  end, let $Q$ be an artinian Gorenstein ring such that $R \is
  Q/\Soc{Q}$. From \eqref{Teter1} one obtains
  \begin{equation}
    \label{eq:IP}
    1 + t\II = (1 + t^3\II)\Poin[Q]{\k} = (1 - t^2 + t^2(1 + t\II))\Poin[Q]{\k}.
  \end{equation}
  Set $1+t\II = \sum_{i=0}^\infty a_it^i$ and $\Poin[Q]{\k} =
  \sum_{i=0}^\infty b_it^i$. Now \eqref{IP} yields relations
  \begin{align*}
    a_0 &= b_0 =1,\\
    a_1 &= b_1,\\
    a_2 &= b_2 -b_0 + b_0a_0 = b_2,
  \end{align*}
  and for $n \ge 2$:
  \begin{equation*}
    a_n = b_n - b_{n-2} + \dsum_{i=0}^{n-2}b_ia_{n-2-i}.
  \end{equation*}
  The coefficients of the series $(1-\ratmin{e}t)(1+t\II)$ may now be
  expressed as
  \begin{align*}
    a_1 - \ratmin{e}a_0 &= b_1 - \ratmin{e}b_0,\\
    a_2 - \ratmin{e}a_1 &= b_2 - \ratmin{e}b_1, \\
    a_3 - \ratmin{e}a_2 &= b_3 - \ratmin{e}b_2 + \ratmin{e}b_0 +
    b_0(a_1 - \ratmin{e}a_0),
  \end{align*}
  and for $n \ge 3$:
  \begin{equation*}
    a_n - \ratmin{e}a_{n-1} = b_n -\ratmin{e}b_{n-1} +
    \ratmin{e}b_{n-3} + \dsum_{i=0}^{n-3}b_i(a_{n-2-i} -
    \ratmin{e}a_{n-3-i}).
  \end{equation*}
  Since $R$ is not a hypersurface, the ring $Q$ is not a hypersurface
  and $\edim{Q} = e$. By \lemref{A} the series $(1 -\ratmin{e}t +
  \ratmin{e}t^3)\Poin[Q]{\k}$ has non-negative coefficients, which
  means that $b_n - \ratmin{e}b_{n-1} + \ratmin{e}b_{n-3} \ge 0$ for
  all $n\ge 1$. It now follows by recursion that $a_n
  -\ratmin{e}a_{n-1} \ge 0$ for all $n\ge 1$.
\end{prf*}


\appendix

\section*{Appendix. Free resolutions over artinian rings}
\label{sec:growth}
\stepcounter{section}

\noindent
Here we collect a few results on the growth of Betti numbers of
modules over an artinian local ring $\Rmk$. Their main application in
this paper is to the module $\Ek$, the injective hull of $\k$.

The first result, \lemref{divides} below, proves the characterization
\pgref{dividesc} of Gorenstein rings. Indeed, let $Q$ be a
Cohen--Macaulay local ring of dimension $d$ and let $\xx =x_1,\ldots,
x_d$ be a $Q$-regular sequence. Then $R=Q/(\xx)$ is artinian, and one
has
\begin{equation*}
  \bass[Q]{d+i} = \bass{i} = \bet{i}{\Ek} \quad\text{for all $i\ge 0$};
\end{equation*}
see \eqref{mb} and \prpcite[1.2.(c)]{HBF71}. If $Q$ is not Gorenstein,
then $R$ is not Gorenstein and \lemref{divides} applies to $\Ek$ and
yields $\bass[Q]{d+i} \ge 2$ for all $i\ge 0$.

\begin{lem}
  \label{lem:divides}
  Let $M$ be a non-free finitely generated $R$-module. If\, $\lgtR$
  divides $\lgt{M}$, then one has $$\bet{i}{M}\ge 2 \ \text{ for all
    $i\ge 0$.}$$ Moreover, if\, $\lgt{M} = \lgtR$, then the inequality
  $\bet{1}{M} \ge \bet{0}{M}$ holds.
\end{lem}

\begin{prf*} 
  Since $M$ is not free, we have $\bet{i}{M}\ge 1$ for all $i\ge 0$.
  For every $j\ge 0$, a length computation based on the exact sequence
  \begin{equation}
    \tag{1}
    \label{eq:star}
    0 \to M_{j+1}\to R^{\bet{j}{M}}\to \cdots\to R^{\bet{0}{M}}\to M\to
    0
  \end{equation}
  shows that $\lgt{M}\equiv\pm\lgt{M_{j+1}}$ (mod $\lgtR$). By
  assumption, $\lgtR$ divides $\lgt{M}$, so it follows that $\lgtR$
  divides $\lgt{M_{j+1}}$; in particular, we have $\lgtR \le
  \lgt{M_{j+1}}$.

  Assume that $\bet{j}{M}=1$ for some $j\ge 0$.  Then the embedding
  $0\to M_{j+1}\to R$ yields $\lgtR > \lgt{M_{j+1}}$, which
  contradicts the inequality obtained above. Therefore, $\bet{i}{M}$
  is at least $2$ for all $i\ge 0$.

  Now assume that $\lgt{M} = \lgtR$. For $j=1$ the sequence
  \eqref{star} gives
  \begin{equation*}
    \lgt{M_2} = (\bet{1}{M}-\bet{0}{M}+1)\lgtR.
  \end{equation*}
  As $\lgt{M_2}>0$, this forces the desired inequality.
\end{prf*}

In the following, $r$ denotes the socle rank of $R$. We set $h =
\max\setof{i}{\m^i \ne 0}$ and adopt the notation from \pgref{art
  notation}. For a finitely generated $R$-module $M$, the rank of the
largest $\k$-vector space that is a direct summand of $M$ is called
the \emph{$\k$-rank} of~$M$.

The next proposition applies to rings of large socle rank, compared to
the length of the ring. It complements a result of Gasharov on Peeva
\prpcite[(2.2)]{VNGIVP90} that applies to rings of large embedding
dimension.

\begin{prp}
  \label{prp:s}
  Assume that $h\ge 2$ and $r > \lgt{R/\m^h}-1$. For every finitely
  generated non-free $R$-module $M$, the sequence $\betseq{M}$ is
  increasing, and for $i\ge 2$ the next inequality holds
  $$\bet{i+1}{M} \ge \frac{r}{\lgt{R/\m^h}-1}\bet{i}{M}.$$
\end{prp}

\begin{prf*} 
  Note that the assumption $h\ge2$ yields $\lgt{R/\m^h}\ge 2$, thus
  the quotient $r/(\lgt{R/\m^h}-1)$ is well-defined.  For $i\ge 1$ set
  $b_i = \bet{i}{M}$ and let $s_{i}$ be the $\k$-rank of
  $M_{i}$. Write the syzygy $M_{i+1}$ as a direct sum $M_{i+1} \is
  \k^{s_{i+1}}\oplus N^{(i+1)}$, where $\k$ is not a direct summand of
  $N^{(i+1)}$. The isomorphism $\SocR^{b_{i}} \is \k^{s_{i+1}}\oplus
  \Soc{N^{(i+1)}}$, see \pgref{soc}, explains the second equality
  below.
  \begin{align*}
    rb_{i} = \rnk{\Soc{R^{b_i}}} &= s_{i+1} + \rnk{\Soc{N^{(i+1)}}} \\
    & \le s_{i+1} + \lgt{\m N^{(i+1)}} \\
    & = s_{i+1} + \lgt{N^{(i+1)}} - \bet{0}{N^{(i+1)}} \\
    &= s_{i+1} + \lgt{N^{(i+1)}} - (b_{i+1} - s_{i+1}) \\
    &\le s_{i+1} + (b_{i+1} - s_{i+1})\lgt[]{R/\m^h} - (b_{i+1} - s_{i+1}) \\
    &= (\lgt[]{R/\m^h} - 1)b_{i+1} + (2- \lgt[]{R/\m^h})s_{i+1}\\
    &\le (\lgt[]{R/\m^h} - 1)b_{i+1}
  \end{align*}
  The first inequality uses the containment $\Soc{N^{(i+1)}} \subseteq
  \m N^{(i+1)}$, which holds as $\k$ is not a direct summand of
  $N^{(i+1)}$; see~\pgref{lescot}. The second inequality follows as
  $N^{(i+1)}$, being a summand of a syzygy, is an $R/\m^h$ module
  generated by $b_{i+1} - s_{i+1}$ elements. The last inequality holds
  as we have $\lgt{R/\m^h}\ge 2$.
\end{prf*}

Applied to the module $\Ek$, the next result establishes the termwise
exponential growth of the sequence $\bassseq[0]$ stated in
\prpref{socgrowth}, see \eqref{mb}.  If there is more than one minimal
generator of the maximal ideal in the socle of $R$, then
\prpref{growth} gives a higher rate of growth than Peeva's
\lemcite[6]{IPv98}.

\begin{prp}
  \label{prp:growth}
  Assume that $R$ satisfies $\m^2 \ne 0$, $\SocR \not\subseteq \m^2$
  and $\edim{R} \ge 2$. For every finitely generated non-free
  $R$-module $M$, the sequence $\betseq[1]{M}$ is increasing, and for
  $i\ge 2$ the next inequality holds
  $$\bet{i+1}{M} \ge (r-\rnk{(\m^2 \cap \SocR)})\left(1+\frac{1}{\lgtR
      -1}\right)\bet{i}{M}.$$
\end{prp}

\begin{prf*}
  For every $i\ge 0$ set $b_i=\bet{i}{M}$. By assumption there is a
  minimal generator of $\m$ in $\SocR$. As $\SocR^{b_1}$ is contained
  in $M_2$, this generator gives $b_1$ elements in $M_2 \setminus \m
  M_2$ that are linearly independent modulo $\m M_2$. It follows from
  the assumptions on $R$ that $\rnk{\SocR}$ is at least $2$, so there
  are elements in $M_2$ which are not in the span of these $b_1$
  elements.  This proves the inequality $b_2 > b_1$.

  Set $A = r-\rnk{(\m^2 \cap \SocR)}$. By assumption $A$ is at least
  $1$, and since $\m^2 \ne 0$ there is a minimal generator of $\m$
  outside the socle, so $e \ge A + 1$.  For $i\ge 2$ write $M_{i}$ as
  a direct sum $M_{i} \is \k^{s_i} \oplus N^{(i)}$, where $\k$ is not
  a summand of $N^{(i)}$.  By \pgref{soc} we have $\Soc{M_i} =
  \SocR^{b_{i-1}}$. A $\k$-vector space of rank $Ab_{i-1}$ is a direct
  summand of $\m R^{b_{i-1}}$ and hence in $\Soc{M_i}\setminus \m
  M_i$. In particular, there is an inequality
  \begin{equation}
    \label{eq:si}
    s_i \ge Ab_{i-1}.
  \end{equation}
  There is an isomorphism $M_{i+1} \is \m^{s_i} \oplus (N^{(i)})_1$,
  so $M_{i+1}$ decomposes as a direct sum $\m^{s_i} \oplus \k^{s}
  \oplus N$, where $\k$ is not a summand of $N$. Notice the equality
  $As_i + s = s_{i+1}$; it explains the second equality in the
  computation below.
  \begin{align*}
    b_{i+1} = es_i + s + \bet{0}{N} &\ge (A+1)s_i + s\\
    &= s_{i+1} + s_i\\
    &\ge Ab_i + Ab_{i-1}\\
    &\ge Ab_i\left(1 + \frac{1}{\lgtR-1}\right)
  \end{align*}
  The penultimate inequality follows from \eqref{si} and the last one
  from \lemcite[4.2.7]{ifr}.
\end{prf*}

\begin{bfhpg}[Rings with radical cube zero]
  \label{JLs}
  Assume $\m^3 = 0 \ne \m^2$ and fix the notation:
  \begin{equation}
    \label{eq:aer}
    a = \rnk{\m^2}, \quad e=\rnk{\m/\m^2}, \qand r=\rnk{\SocR}.
  \end{equation}
  \noindent
  Notice that the inclusion $\m^2\subseteq\SocR$ yields the inequality
  $a\le r$.

  Let $M$ be a finitely generated $R$-module, and set $b_i =
  \bet{i}{M}$ for $i\ge 0$. We recall a few facts from
  Lescot's~\cite[\S 3]{JLs85}. There is an inequality
  \begin{equation}
    \label{eq:betti1}
    b_{1} \ge eb_{0} - \rnk{(\m M/\m^2M)}.
  \end{equation}
  If $\k$ is not a direct summand of $M_1$, then equality holds; the
  converse is true if $\SocR = \m^2$.

  For every $i\ge 2$ there is an inequality
  \begin{equation}
    \label{eq:betti3} b_{i} \ge eb_{i-1} - ab_{i-2}.
  \end{equation}
  
  If $\k$ is not a summand of $M_1$ and not a summand of $M_2$, then
  the equality
  \begin{equation}
    b_{2} = eb_{1} - ab_{0}
  \end{equation}
  holds if and only if $\m^2M=0$.

  If $i\ge 3$, and $\k$ is not a summand of $M_i$ and not a summand of
  $M_{i-1}$, then
  \begin{equation}
    \label{eq:betti eq}
    b_{i} = eb_{i-1} - ab_{i-2}.
  \end{equation}

  If $\SocR = \m^2$ and $\m^2M=0$, then
  \begin{equation}
    \label{eq:betti2}
    eb_{1}\ge  rb_{0} + s(e-1),
  \end{equation}
  where $s$ is the $\k$-rank of $M_1$.
\end{bfhpg}

The next result strengthens parts of Lescot's \thmcite[B]{JLs85}.

\begin{thm}
  \label{thm:increasing}
  Let $R$ be a local ring with $\m^3=0\not=\m^2$ that is not
  Gorenstein. Let $M$ be a finitely generated $R$-module that
  satisfies $\bet{1}{M}>\bet{0}{M}$. With the notation from {\rm
    \eqref{aer}} the following statements hold
  \begin{prt}
  \item If $a\ne e$ or $a\ne r$, then the sequence $\betseq[0]{M}$ is
    increasing and has termwise exponential growth.
  \item If $a=e=r$, then the sequence $\betseq[0]{M}$ is
    non-decreasing with strict inequalities $\bet{i}{M} >
    \bet{i-1}{M}$ for all but, possibly, one index $i$. Moreover, if
    there is an equality $\bet{j}{M} = \bet{j-1}{M}$, then $j$ is at
    least $2$, and $\k$ is a direct summand of $M_{j+1}$ and not a
    direct summand of $M_i$ when $0\le i \le j$.
  \end{prt}
\end{thm}

\begin{prf*}
  For $i\ge 0$ set $b_i = \bet{i}{M}$.

  (a): First assume $a\ne r$, that is, $\m^2\not=\SocR$. By
  \prpcite[3.9]{JLs85} and the assumption on $M$, the sequence $\bseq$
  is increasing, and by \prpref{growth} it has termwise exponential
  growth of rate $(r-a)(1 + 1/(a+e))$.

  Next, assume $a=r$. There are two cases to consider:

  Case $e<a$: Since $\lgt{R/\m^2}-1=e$, it follows from \prpref{s} and
  the assumption on $M$ that the sequence $\bseq$ is increasing and
  has termwise exponential growth of rate $a/e$.
  
  Case $e>a$: For every $i\ge 2$, the inequality \eqref{betti3} yields
  \begin{equation*}
    b_{i} -b_{i-1} \ge (e-1)b_{i-1} - ab_{i-2} \ge
    a(b_{i-1}-b_{i-2}). 
  \end{equation*}  
  By recursion, based on the assumption $b_1 > b_0$, it follows that
  the sequence $\bseq$ is increasing. If $e=a+1$, then the sequence
  has termwise exponential growth by \prpcite[3]{IPv98}. If $e > a+1$,
  then the same conclusion follows as \eqref{betti3} yields
  \begin{equation*}
    b_{i} > (e-a)b_{i-1} \quad\text{for all $i\ge 2$}. %
  \end{equation*} 
  
  (b): Assume that $a=e=r$.  First notice that if $i \ge 3$ and $\k$
  is not a summand of $M_i$ and not a summand of $M_{i-1}$, then
  \eqref{betti eq} yields $b_{i-1} > b_{i-2}$, as $b_i>0$ by
  assumption. It follows that the sequence $\bseq$ is increasing if
  for every $i \ge 2$ the residue field $\k$ is not a direct summand
  of $M_i$.

  Next, assume that $\k$ is a summand of a syzygy of $M$, and let $j$
  be the least integer such that $M_{j+1}=\k\oplus N$ for some
  $R$-module $N$.  The sequence $\bseq[j+1]$ is then
  increasing. Indeed, the sequence $\betseq[0]{\k}$ is increasing, as
  $R$ is not a hypersurface. If $j\ge 0$, then the sequence
  $\betseq[0]{N}$ is non-decreasing by \eqref{betti2}, as $\m^2
  N=0$. If $j=-1$, i.e.\ $\k$ is a summand of $M$, then the sequence
  $\betseq[1]{N}$ is non-decreasing by \eqref{betti2}, and $b_1 > b_0$
  by assumption.  If $j=-1$ or $j=0$, it is thus immediate that the
  sequence $\bseq$ is increasing.

  We can now assume that $j$ is at least $1$. In the next chain of
  inequalities, the right-most and left-most ones are already know;
  the inequalities in-between follow by application of \eqref{betti2}
  to the syzygies $M_1, \dots, M_j$.
  \begin{equation*}
    \cdots > b_{j+1} > b_j \ge \cdots \ge b_{1} > b_0
  \end{equation*}
  If $j=1$, it follows that the sequence $\bseq$ is increasing.  If
  $j\ge 2$, then $b_{i-1} > b_{i-2}$ for $i$ with $j \ge i \ge 2$ by
  \eqref{betti eq} and the assumption on $M$. In total, this gives
  inequalities $\cdots > b_{j+1} > b_{j} \ge b_{j-1} > \cdots > b_{1}
  > b_{0}$.\qedhere
\end{prf*}

\section*{Acknowledgments}

We thank Luchezar L.\ Avramov and Sean Sather-Wagstaff for extensive
comments on an earlier version of the paper. Thanks are also due to
Christopher Monico and Roger W.~Barnard for discussions related to the
material in \secref{teter}. A series of pertinent remarks from the
anonymous referee helped us improve the exposition.

\enlargethispage*{2\baselineskip}
\bibliographystyle{amsplain}

\newcommand{\arxiv}[2][AC]{\mbox{\href{http://arxiv.org/abs/#2}{\sf
      arXiv:#2 [math.#1]}}}
\newcommand{\oldarxiv}[2][AC]{\mbox{\href{http://arxiv.org/abs/math/#2}{\sf
      arXiv:math/#2
      [math.#1]}}}\providecommand{\MR}[1]{\mbox{\href{http://www.ams.org/mathscine%
      t-getitem?mr=#1}{#1}}}
\renewcommand{\MR}[1]{\mbox{\href{http://www.ams.org/mathscinet-getitem?mr=#%
      1}{#1}}} \providecommand{\bysame}{\leavevmode\hbox
  to3em{\hrulefill}\thinspace}
\providecommand{\MR}{\relax\ifhmode\unskip\space\fi MR }
\providecommand{\MRhref}[2]{%
  \href{http://www.ams.org/mathscinet-getitem?mr=#1}{#2} }
\providecommand{\href}[2]{#2}

\end{document}